\documentclass[12pt]{amsart}
\usepackage{amssymb,amsmath,amsthm,mathtools,amsfonts,times,hyperref,mathrsfs,multirow}

\newtheorem{theorem}{Theorem}[section]

\theoremstyle{definition}
\newtheorem{definition}[theorem]{Definition}
\newtheorem{example}[theorem]{Example}

\theoremstyle{remark}
\newtheorem{remark}{Remark}
\numberwithin{equation}{section}

\allowdisplaybreaks

\newcommand\thm[1]{\ref{thm:#1}}

\newcommand\eqn[1]{(\ref{eq:#1})}

\newcommand\refdef[1]{\ref{def:#1}}

\begin{document}

\title[\scalebox{.9}{mex of overpartitions}]{On new minimal excludants of overpartitions related to some $q$-series of Ramanujan}

\author{Aritram Dhar}
\address{Department of Mathematics, University of Florida, Gainesville
FL 32611, USA}
\email{aritramdhar@ufl.edu}
\author{Avi Mukhopadhyay}
\address{Department of Mathematics, University of Florida, Gainesville
FL 32611, USA}
\email{mukhopadhyay.avi@ufl.edu}
\author{Rishabh Sarma}
\address{Department of Mathematics, University of Florida, Gainesville
FL 32611, USA}
\email{rishabh.sarma@ufl.edu}

\date{\today}



\begin{abstract}
Inspired by Andrews' and Newman's work on the minimal excludant or ``mex'' of partitions, we define four new classes of minimal excludants for overpartitions and establish relations to certain functions due to Ramanujan.
\end{abstract}
\maketitle

\section{Introduction}
\label{section1}

Given an integer partition $\pi$, the minimal excludant of $\pi$ is defined to be the smallest positive integer that is not a part of $\pi$. This partition statistic seems to be first considered by Grabner and Knopfmacher \cite{Gr-Kno06} in 2006, who call it the least gap. They obtained the result that the sum of the minimal excludants over all partitions of $n$ is equal to is equal to the number of partitions of $n$ into distinct parts with two colors. In 2019, Andrews and Newman introduced the terminology ``minimal excludant'' or ``mex'' of an integer partition function and initiated the study of the connections of the mex to other partition theoretic objects and statistics. Among these results is a rediscovery of the above result of Grabner and Knopfmacher. This statistic has also appeared earlier in 2011 in a work of Andrews \cite{An11} where he relates the minimal excludant statistic, then called the smallest part of a partition that is not a summand by Andrews, to the Frobenius symbol representation of partitions.
\\\\
In the first of a series of two papers, Andrews and Newman \cite{An-New19} considered $\sigma$\,mex($n$), the sum of mex($\pi$) taken over all partitions $\pi$ of $n$. They showed that $\sigma$\,mex($n$) is equal to the number of partitions of $n$ into distinct parts with two colors \cite[Theorem 1.1]{An-New19}. Analogously, Aricheta and Donato \cite{Ar-Don23} have extended this concept to overpartitions. To recall, an overpartition is a partition in which the first occurrence of a number may be overlined \cite{Co-Lov04}. Aricheta and Donato define the minimal excludant or mex of an overpartition $\pi$, denoted by $\overline{\text{mex}}(\pi)$, to be the smallest positive integer that is not a part of the nonoverlined parts of $\pi$. Following this, for a positive integer $n$, the authors define $\sigma$\,$\overline{\text{mex}}(n)$ to be the sum of the minimal excludants over all overpartitions $\pi$ of $n$ and prove that $\sigma$\,$\overline{\text{mex}}(n)$ equals the number of partitions of $n$ into distinct parts using three colors \cite[Theorem 1.1]{Ar-Don23}.
\\\\
In the second paper on the minimal excludant statistic, Andrews and Newman \cite{An-New20} defined an extended function of their minimal excludant of a partition defined above, and explored relations of this extended mex function to other well studied partition statistics such as the rank and crank. Furthermore, the connection between the mex and the crank made by Andrews and Newman was independently made by Hopkins and Sellers \cite{Ho-Sel20} in the same year. Subsequently, generalized versions of these relations were studied by us \cite{Dh-Mu-Sar23} by calculating the generating function of the general case of Andrews' and Newman's extended minimal excludant. At the end of the very same paper \cite[Section 4]{Dh-Mu-Sar23}, we had defined a new minimal excludant for overpartitions which to our surprise was related to a function of Ramanujan. In this paper, we continue our study on minimal excludants of overpartitions and their relations to two fifth order mock theta functions and some other $q$-series related to Ramanujan.
\\\\
Let $L,m,n$ be non-negative integers. Throughout the paper, we use the following standard notations \cite{An98}.
\begin{align*}
    (a)_L = (a;q)_L &:= \prod_{k=0}^{L-1}(1-aq^k),\\
    (a)_{\infty} = (a;q)_{\infty} &:= \lim_{L\rightarrow \infty}(a)_L\,\,\text{where}\,\,\lvert q\rvert<1.
\end{align*}
We define the $q$-binomial coefficient as
\begin{align*}
    \left[\begin{matrix}m\\n\end{matrix}\right]_q :=  \Bigg\{\begin{array}{lr}
        \dfrac{(q)_m}{(q)_n(q)_{m-n}}\quad\text{for } m\ge n\ge 0,\\
        0\qquad\qquad\quad\text{otherwise}.\end{array}
\end{align*}
\\
Let $\overline{P}(q) = \dfrac{(-q)_{\infty}}{(q)_{\infty}}$ and $\overline{P}_\text{o}(q) = \dfrac{(-q;q^2)_{\infty}}{(q;q^2)_{\infty}}$ denote the generating functions of overpartitions and overpartitions into odd parts respectively.
\\\\
We define four new classes of minimal excludants of overpartitions with an example for each. The first of these minimal excludants appears in our previous paper \cite[Section 4]{Dh-Mu-Sar23}.
\subsection{Four overpartition mex statistics}
\begin{definition}
\label{def:mex1}
For an overpartition $\pi$, let omex$(\pi)$ be the smallest positive integer that is not a part (overlined or non-overlined) of $\pi$.
\end{definition}
\begin{remark}
One can observe that the minimal excludant in the above definition differs from the one considered by Aricheta and Donato in the light of the fact that their mex is the smallest part missing from the non-overlined parts of the partition. Whereas in our case, we are free to pick the smallest missing part from either the overlined or non-overlined parts of the partition in consideration.
\end{remark}
\begin{example}
If we consider the following two overpartitions
\begin{align*}
\pi_1&=\overline{5}+\overline{4}+4+2+1,
\\
\pi_2&=\overline{10}+8+5+\overline{3}+2+\overline{1},
\end{align*}
then, omex$(\pi_1)=3$ and omex$(\pi_2)=4$.
\end{example}
\noindent
\textbf{A $q$-series generating function related to omex$(\pi)$} : Let $\overline{m}(n)$ denote the number of overpartitions $\pi$ of $n$ having the property that no positive integer less than  omex$(\pi)$ is overlined. As an example, below are the eight overpartitions of $3$ :
\begin{align*}
3, \overline{3}, 2+1, \overline{2}+1, 2+\overline{1}, \overline{2}+\overline{1}, 1+1+1, \overline{1}+1+1.
\end{align*}
\noindent
Then, the overpartitions $\pi \in \{3, \overline{3}, 2+1, 1+1+1\}$ satisfy the stated property, and hence $\overline{m}(3)=4$. Let $$\overline{M}(q)=\sum\limits_{n=0}^{\infty}\overline{m}(n)q^n.$$

\begin{definition}
\label{def:mex2}
For an overpartition $\pi$ into odd parts, let omoex$(\pi)$ be the smallest positive odd integer which is not a part (overlined or non-overlined) of $\pi$.
\end{definition}

\begin{example}
For the overpartitions
\begin{align*}
\pi_3&=\overline{7}+\overline{7}+3+1,
\\
\pi_4&=\overline{7}+7+\overline{5}+3+1,
\end{align*}
we have, omoex$(\pi_3)=5$ and omoex$(\pi_4)=9$.
\end{example}
\noindent
\textbf{A $q$-series generating function related to omoex$(\pi)$} : Let $\overline{m}_o(n)$ denote the number of overpartitions $\pi$ of $n$ into odd parts having the property that no positive integer less than omoex$(\pi)$ is overlined. For example, for $n=3$, the overpartitions $\pi \in \{3, \overline{3}, 1+1+1\}$ satisfy the stated property, and hence $\overline{m}_o(3)=3$. Let $$\overline{M}_\text{o}(q)=\sum\limits_{n=0}^{\infty}\overline{m}_o(n)q^n.$$
\\
For our other two classes of overpartition mexes, we define an ordering on the parts of an overpartition $\pi$, where every non-overlined part is smaller than its overlined counter-part i.e.
$$1<\overline{1}<2<\overline{2}<3<\overline{3}<\cdots$$
\noindent
In the two definitions that follow, we use a $\sim$ symbol over our respective mexes to denote the minimal excludant taken over overpartitions with the ordering on the parts defined above.  

\begin{definition}
\label{def:mexord1}
For an overpartition $\pi$ where we take into consideration the aforementioned ordering on the parts, let $\widetilde{\text{omex}}(\pi)$ be the smallest overlined positive integer which is not a part of $\pi$.
\end{definition}

\begin{example}
For the overpartitions
\begin{align*}
\pi_5&=\overline{5}+\overline{3}+2+1,
\\
\pi_6&=\overline{7}+7+\overline{5}+3+\overline{2}+\overline{1},
\end{align*}
we have, $\widetilde{\text{omex}}(\pi_5)=\overline{1}$ and $\widetilde{\text{omex}}(\pi_6)=\overline{3}$.
\end{example}
\noindent
\textbf{A $q$-series generating function related to $\widetilde{\text{omex}}(\pi)$} : 
Let $\widetilde{m}(n)$ denote the number of overpartitions $\pi$ of $n$ having the property that all overlined and non-overlined parts smaller than  $\widetilde{\text{omex}}(\pi)$ occur as parts. For example, for $n=3$, the overpartitions $\pi \in \{2+1,\overline{2}+1 ,1+1+1\}$ satisfy the stated property, and hence $\widetilde{m}(3)=3$. Let $$\widetilde{M}(q)=\sum\limits_{n=0}^{\infty}\widetilde{m}(n)q^n.$$

\begin{definition}
\label{def:mexord2}
For an overpartition $\pi$ into odd parts where we once again take into consideration the same ordering on the parts as before, let $\widetilde{\text{omoex}}(\pi)$ be the smallest overlined positive odd integer which is not a part of $\pi$.
\end{definition}

\begin{example}
For the overpartitions
\begin{align*}
\pi_7&=\overline{5}+\overline{3}+3+\overline{1},
\\
\pi_8&=11+\overline{7}+\overline{5}+\overline{3}+3+\overline{1},
\end{align*}
we have, $\widetilde{\text{omoex}}(\pi_7)=\overline{7}$ and $\widetilde{\text{omoex}}(\pi_8)=\overline{9}$.
\end{example}
\noindent
\textbf{A $q$-series generating function related to $\widetilde{\text{omoex}}(\pi)$} : 
Let $\widetilde{m}_o(n)$ denote the number of overpartitions of $n$ into odd parts having the property that all overlined and non-overlined parts smaller than  $\widetilde{\text{omoex}}(\pi)$ occur as parts. For example, the only overpartition of $3$ that satisfies the stated property is $1+1+1$, and hence $\widetilde{m}_o(3)=1$. Let $$\widetilde{M}_\text{o}(q)=\sum\limits_{n=0}^{\infty}\widetilde{m}_o(n)q^n.$$

\subsection{Two arithmetic mex functions}
Finally, analogous to the arithmetic function $$\sigma {\text{mex}}(n) = \sum_{\pi\vdash n} {\text{mex}}(\pi)$$ over partitions $\pi$ of $n$ considered by Andrews and Newman, we study the following two analogous sums of our first two minimal excludants of overpartitions introduced in Definitions \refdef{mex1} and \refdef{mex2}. To that end, let us consider the sum 
$$\sigma{\text{omex}}(n) = \sum_{\pi\vdash n} {\text{omex}}(\pi)$$ 
taken over all overpartitions $\pi$ of $n$. Define $\overline{M}(z,q)$ to be the double series in which the coefficient of $z^mq^n$ is the number of overpartitions $\pi$ of $n$ with ${\text{omex}}(\pi) = m$ and let 
$$\sigma \overline{M}(q) = \sum_{n\ge 0}\sigma {\text{omex}}(n)q^n.$$
Again, we consider the sum 
$$\sigma {\text{omoex}}(n) = \sum_{\pi\vdash n} {\text{omoex}}(\pi)$$
taken over all overpartitions $\pi$ of $n$ into odd parts. Define $\overline{M}_o(z,q)$ to be the double series in which the coefficient of $z^mq^n$ is the number of overpartitions $\pi$ of $n$ into odd parts with ${\text{omoex}}(\pi) = m$ and let 
$$\sigma \overline{M}_\text{o}(q) = \sum_{n\ge 0}\sigma {\text{omoex}}(n)q^n.$$

\section{Main Results}
\label{section2}
In this section, we present the statements of our results. Theorem \thm{mex1} below is also stated and proved in our paper on the generalization of the mex of partitions statistic (see \cite[Theorem 30]{Dh-Mu-Sar23} for details). We have included the statement of the result here for the sake of completion.
\begin{theorem}
\label{thm:mex1}
\begin{align*}
\overline{M}(q) = \overline{P}(q)(2 - R(q)),   
\end{align*}
where 
$$R(q) = {\displaystyle\sum\limits_{n=0}^{\infty}\dfrac{q^{\frac{n(n+1)}{2}}}{(-q)_n}} = 1 + \sum\limits_{n=1}^{\infty}(-1)^{n-1}q^n(q)_{n-1} = (q)_{\infty} + 2(q)_{\infty}\sum\limits_{n=1}^{\infty}\dfrac{q^n}{(q)_n(1+q^n)}.$$
\end{theorem}

\begin{remark}
The first two $q$-series representations of $R(q)$ above are due to Ramanujan \cite{Ra88} and the last representation is due to Gupta \cite[Equation 1.11]{Gu21}. It was studied by Andrews \cite{An86} in connection with identities from Ramanujan's ``Lost'' Notebook \cite{Ra88} that involve $R(q)$, the generating functions of the number of divisors of $n$ and that of the number of partitions of $n$ into distinct parts. Subsequently, conjectures made by Andrews in this paper on the distribution of the coefficients of $R(q)$ were proved by Andrews, Dyson and Hickerson \cite{An-Dy-Hic88}. 
\end{remark}

\begin{theorem}
\label{thm:mex2}
\begin{align*}
\overline{M}_\text{o}(q) = \overline{P}_\text{o}(q)(1 - F(-q)),    
\end{align*}
where 
$$F(q) = {\displaystyle\sum\limits_{n=1}^{\infty}\dfrac{(-1)^n{q^{n^2}}}{(q;q^2)_n}}$$ 
is a companion function to $R(q)$.
\end{theorem}

\begin{remark}
Consider partitions into odd parts, with the property that if $k$ occurs as a part, then all positive odd numbers less than $k$ also occur. Then $F(q)$ is the generating function for the number of such partitions where the largest part is congruent to $3$ modulo $4$ minus the number of such partitions where the largest part is congruent to $1$ modulo $4$. See \cite[Section 5]{An-Dy-Hic88} for a treatment of $F(q)$. 
\end{remark}

\begin{theorem}
\label{thm:mexord1}
\begin{align*}
\widetilde{M}(q) = \overline{P}(q)(f_0(q) - 1),    
\end{align*}
where 
$$f_0(q) = {\displaystyle\sum\limits_{n=0}^{\infty}\dfrac{q^{n^2}}{(-q)_n}}$$ 
is a fifth order mock theta function of Ramanujan.
\end{theorem}

\begin{theorem}
\label{thm:mexord2}
\begin{align*}
\widetilde{M}_\text{o}(q) = q\overline{P}_\text{o}(q)F_1(-q),    
\end{align*}
where 
$$F_1(q) = {\displaystyle\sum\limits_{n=0}^{\infty}\dfrac{q^{2n^2+2n}}{(q;q^2)_{n+1}}}$$ 
is a fifth order mock theta function of Ramanujan.
\end{theorem}

\begin{remark}
Ramanujan’s mock theta functions are examples of mock modular forms. Functions of $\widetilde{M}(q)$ and $\widetilde{M}_\text{o}(q)$ type that arise from by multiplying or dividing by certain infinite products are known as mixed mock modular forms and possess interesting properties. See \cite{Lo-Os13} for a treatment of mixed mock modular forms.
\end{remark}

In the following theorem, we relate our sum of mex function $\sigma\overline{M}(q)$ defined in the introduction section to a $q$-series where a finite number of starting consecutive terms are subtracted from an infinite $q$-series to obtain a sum-of-tails series.

\begin{theorem}
\label{thm:mexsum1}
\begin{align*}
\sigma\overline{M}(q) = \overline{P}(q)\left(R(q)-2(q)_{\infty}\sum\limits_{n=0}^{\infty}\dfrac{q^n}{(q)_n(1+q^n)}G_n(q)\right),   
\end{align*}
where $G_n(q)$ is the $q$-series tail given by 
$$G_n(q)=\sum\limits_{i=n+1}^{\infty}\dfrac{q^i}{1+q^i}.$$ 
\end{theorem}

\begin{theorem}
\label{thm:mexsum2}
\begin{align*}
\sigma\overline{M}_\text{o}(q) = \overline{P}_\text{o}(q)\left(1+q\sum_{n=1}^{\infty}(-1)^n(q^2;q^2)_nq^nH_n(q^2)\right),     
\end{align*}
where the partial sum 
$$H_n(q)=\sum_{i=1}^{n}\frac{q^{i}}{1-q^{i}}$$ 
is the $q$-analog of the harmonic number $H_n=\sum\limits_{i=1}^{n}\frac{1}{i}$.
\end{theorem}

\begin{remark}
The $q$-harmonic series $H_n(q)$ above are indeed partial sums of the divisor function generating function given by 
$$\sum_{i=1}^{\infty}d(i)q^i=\sum_{i=1}^{\infty}\frac{q^{i}}{1-q^{i}}=\sum_{i=1}^{\infty}\frac{(-1)^{i-1}q^{\frac{i(i+1)}{2}}}{{(1-q^i)(q;q)_i}}.$$
Interesting formulas for harmonic and $q$-harmonic numbers $H_n$ and $H_n(q)$ were re-established by Andrews and Uchimura in \cite{An-Uch85} using differentiation of classical hypergeometric series technique. One such formula relevant to this discussion is the finite analog of the generating function of the divisor function viz. $$H_n(q)=\sum_{i=1}^{n}\frac{q^{i}}{1-q^{i}}=\sum_{i=1}^{n}\frac{(-1)^{i-1}}{1-q^i}q^{\frac{i(i+1)}{2}}\left[\begin{matrix}n\\i\end{matrix}\right]_q.$$ 
\end{remark}

\section{Proofs of Main Results}
\label{section3}
In this section, we provide proofs of our results. We follow the notation of Gasper and Rahman \cite{Ga-Rah04}. For instance, the ${}_r\phi_{r-1}$ unilateral basic hypergeometric series with base $q$ and argument $z$ is defined by
\begin{align*}
{}_r\phi_{r-1} \left(
	\setlength\arraycolsep{2pt}
	\begin{matrix}
		a_1,\ldots,a_r \\
		\multicolumn{2}{c}{
			\begin{matrix}
				b_1,\ldots,b_{r-1} 	
			\end{matrix}}
	\end{matrix} \hspace{1pt}
;q, z \right) &:= \sum_{k=0}^{\infty}\dfrac{(a_1,\ldots,a_r;q)_k}{(q,b_1,\ldots,b_{r-1};q)_k}z^k,\quad \lvert z\rvert<1. 
\end{align*}
We also need the $q$-binomial theorem : For $|z|<1$, 
\begin{align}
\label{eq:qbinom}
{}_1\phi_0\left( \begin{gathered} a \\ - \end{gathered} ;\,q, z \right)=\sum\limits_{k=0}^{\infty}\frac{(a;q)_k}{(q;q)_k}z^k=\frac{(az;q)_{\infty}}{(z;q)_{\infty}},
\end{align}
and Heine's ${}_2\phi_1$ transformation : For $|z|<1$ and $|b|<1$, 
\begin{align}
\label{eq:heine}
{}_2\phi_1\left( \begin{gathered} a, b \\ c \end{gathered}
;\,q, z \right)=\frac{(b;q)_{\infty}(az;q)_{\infty}}{(c;q)_{\infty}(z;q)_{\infty}}{}_2\phi_1\left( \begin{gathered} c/b, z \\ az \end{gathered}
;\,q, b \right).
\end{align}
Finally, we will also make use of the following simple identity \cite[Equation 3]{Dh-Mu-Sar23} that we had obtained from a ${}_1\phi_1$ summation of Gasper-Rahman and used to prove Theorem \thm{mex1} in \cite[Section 4]{Dh-Mu-Sar23}. 
\begin{align}
\label{eq:gasrah}
\sum\limits_{n=0}^{\infty}\frac{z^nq^{\frac{n(n-1)}{2}}}{(-zq;q)_n}=1+z.
\end{align}
\subsection{Proof of Theorem \thm{mex2}}
By standard combinatorial arguments, we deduce that
\begin{align*}
\overline{M}_\text{o}(q) &= \sum\limits_{n=0}^{\infty}\overline{m}_o(n)q^n
\\
&= \sum\limits_{n=1}^{\infty}\dfrac{q^{1+3+\cdots+(2n-3)}\prod\limits_{m=n+1}^{\infty}(1+q^{2m-1})}{\prod\limits_{\substack{m=1 \\ m \neq n}}^{\infty}(1-q^{2m-1})}
\\
&= \frac{(-q;q^2)_{\infty}}{(q;q^2)_{\infty}}\sum\limits_{n=1}^{\infty}\dfrac{q^{(n-1)^2}(1-q^{2n-1})}{(-q;q^2)_n}
\\
&= \frac{(-q;q^2)_{\infty}}{(q;q^2)_{\infty}}\left[\sum\limits_{n=1}^{\infty}\dfrac{q^{(n-1)^2}}{(-q;q^2)_n}-\sum\limits_{n=1}^{\infty}\dfrac{q^{n^2}}{(-q;q^2)_n}\right]
\\
&= \overline{P}_\text{o}(q)\left(1-F(-q)\right)
\end{align*}
where the last line follows by replacing $q\mapsto q^2$ and substituting $z = q$ in Equation \eqn{gasrah}.\qed

\subsection{Proof of Theorem \thm{mexord1}}
By standard combinatorial arguments, we deduce that
\begin{align*}
\widetilde{M}(q) &= \sum\limits_{n=1}^{\infty}\widetilde{m}(n)q^n
\\
&= \sum\limits_{n=1}^{\infty}\dfrac{q^{1+1+2+2+\cdots+(n-1)+(n-1)+n}\prod\limits_{m=n+1}^{\infty}(1+q^{m})}{\prod\limits_{\substack{m=1}}^{\infty}(1-q^{m})}
\\
&= \frac{(-q)_{\infty}}{(q)_{\infty}}\sum\limits_{n=1}^{\infty}\dfrac{q^{n^2}}{(-q)_n}
\\
&= \overline{P}(q)\left[f_0(q)-1\right].
\end{align*}\qed

\subsection{Proof of Theorem \thm{mexord2}}
By standard combinatorial arguments, we deduce that
\begin{align*}
\widetilde{M}_\text{o}(q) &= \sum\limits_{n=1}^{\infty}\widetilde{m}_o(n)q^n
\\
&= \sum\limits_{n=1}^{\infty}\dfrac{q^{1+1+3+3+\cdots+(2n-3)+(2n-3)+(2n-1)}\prod\limits_{m=n+1}^{\infty}(1+q^{2m-1})}{\prod\limits_{\substack{m=1}}^{\infty}(1-q^{2m-1})}
\\
&= \frac{(-q;q^2)_{\infty}}{(q;q^2)_{\infty}}\sum\limits_{n=1}^{\infty}\dfrac{q^{2n^2-2n+1}}{(-q;q^2)_n}
\\
&= \frac{(-q;q^2)_{\infty}}{(q;q^2)_{\infty}}\sum\limits_{n=0}^{\infty}\dfrac{q^{2n^2+2n+1}}{(-q;q^2)_{n+1}}
\\
&= q\overline{P}_\text{o}(q)F_1(-q).
\end{align*}\qed

\subsection{Proof of Theorem \thm{mexsum1}}
The proof follows from standard combinatorial arguments and using Heine's ${}_2\phi_1$ transformation stated at the beginning of the section. During the course of our $q$-series transformations, we also need the following identity due to Gupta \cite[Equation 1.15]{Gu21} valid for any $c \in \mathbb{R}$ and $|t|<1$.
\begin{align}
\label{eq:rajat}
\sum\limits_{n=0}c^n\left((t)_n - (t)_{\infty}\right) = (t)_{\infty}\sum\limits_{n=1}^{\infty}\dfrac{t^n}{(q)_n(1-cq^n)}.
\end{align}
\noindent
Proceeding with the proof of our theorem, we have
\begin{align*}
\overline{M}(z,q) &= \sum\limits_{n=1}^{\infty}\dfrac{z^nq^{1+2+\cdots+(n-1)}\prod\limits_{m=n+1}^{\infty}(1+q^{m})}{\prod\limits_{\substack{m=1 \\ m \neq n}}^{\infty}(1-q^{m})}
\\
&= \frac{(-q)_{\infty}}{(q)_{\infty}}\sum\limits_{n=1}^{\infty}\dfrac{z^nq^{\frac{n(n-1)}{2}}(1-q^{n})}{(-q)_n}
\\
&= \overline{P}(q)\left[\sum\limits_{n=0}^{\infty}\dfrac{z^nq^{\frac{n(n-1)}{2}}}{(-q)_n}-\sum\limits_{n=0}^{\infty}\dfrac{z^nq^{\frac{n(n+1)}{2}}}{(-q)_n}\right].
\end{align*}
Thus
\begin{align*}
\sigma\overline{M}(q) &= \dfrac{\partial}{\partial z}\Big\vert_{z=1}\overline{M}(z,q)
\\
&= \overline{P}(q)\left[A(q)-B(q)\right],
\end{align*}
where
\begin{align*}
A(q) &= \dfrac{\partial}{\partial z}\Big\vert_{z=1}\sum\limits_{n=0}^{\infty}\dfrac{z^nq^{\frac{n(n-1)}{2}}}{(-q)_n}
\\
&= \dfrac{\partial}{\partial z}\Big\vert_{z=1}\lim\limits_{\tau\rightarrow 0}{}_2\phi_1\left( \begin{gathered} -1/\tau, q \\ -q \end{gathered}
;\,q, z \tau \right)
\\
&= \dfrac{\partial}{\partial z}\Big\vert_{z=1}\lim\limits_{\tau\rightarrow 0}\dfrac{(q)_{\infty}(-z)_{\infty}}{(-q)_{\infty}(z\tau)_{\infty}}{}
_2\phi_1\left( \begin{gathered} -1, z\tau \\ -z \end{gathered}
;\,q, q \right)
\\
&\quad\text{(using Equation \eqn{heine} with $(a, b, c, z) \mapsto (-1/\tau, q, -q, z\tau)$)}
\\
&= \dfrac{\partial}{\partial z}\Big\vert_{z=1}\dfrac{(q)_{\infty}}{(-q)_{\infty}}\sum\limits_{n=0}^{\infty}\dfrac{(-1)_n(-zq^n)_{\infty}q^n}{(q)_n}
\\
&= 2(q)_{\infty}\sum\limits_{n=0}^{\infty}\dfrac{q^n}{(q)_n}\sum\limits_{i=n}^{\infty}\dfrac{q^i}{1+q^i}
\\
&= 2(q)_{\infty}\sum\limits_{n=1}^{\infty}\dfrac{q^n}{(q)_n}\sum\limits_{i=n}^{\infty}\dfrac{q^i(-q)_{i-1}}{(-q)_i}+2(q)_{\infty}\sum\limits_{i=0}^{\infty}\dfrac{q^i}{1+q^i}
\\
&= 2(q)_{\infty}\sum\limits_{n=1}^{\infty}\dfrac{q^{2n}}{(q)_n}\sum\limits_{i=0}^{\infty}\dfrac{q^i(-q)_{i+n-1}}{(-q)_{i+n}}+2(q)_{\infty}\sum\limits_{i=0}^{\infty}\dfrac{q^i}{1+q^i}
\\
&= 2(q)_{\infty}\sum\limits_{n=1}^{\infty}\dfrac{q^{2n}(-q)_{n-1}}{(q)_n(-q)_n}{}
_2\phi_1\left( \begin{gathered} q,-q^n \\ -q^{n+1} \end{gathered}
;\,q, q \right)
+2(q)_{\infty}\sum\limits_{i=0}^{\infty}\dfrac{q^i}{1+q^i}
\\
&= 2(q^2)_{\infty}\sum\limits_{n=1}^{\infty}\dfrac{q^{2n}}{(q)_n}\sum\limits_{i=0}^{\infty}\dfrac{(-1)^i(q)_iq^{ni}}{(q^2)_i}+2(q)_{\infty}\sum\limits_{i=0}^{\infty}\dfrac{q^i}{1+q^i}
\\
&\quad\text{(using Equation \eqn{heine} with $(a, b, c, z) \mapsto (q, -q^n, -q^{n+1}, q)$ in the first double sum)}
\\
&= 2(q^2)_{\infty}\sum\limits_{i=0}^{\infty}\dfrac{(-1)^i(q)_i}{(q^2)_i}\sum\limits_{n=1}^{\infty}\dfrac{q^{(i+2)n}}{(q)_n}+2(q)_{\infty}\sum\limits_{i=0}^{\infty}\dfrac{q^i}{1+q^i}
\\
&= 2(q^2)_{\infty}\sum\limits_{i=0}^{\infty}\dfrac{(-1)^i(q)_i}{(q^2)_i}\left[\dfrac{1}{(q^{i+2})_{\infty}}-1\right]+2(q)_{\infty}\sum\limits_{i=0}^{\infty}\dfrac{q^i}{1+q^i}
\\
&\quad\text{(using Equation \eqn{qbinom} with $(a, z) \mapsto (0, q^{i+2})$ in the first double sum)}
\\
&= 2\sum\limits_{n=0}^{\infty}(-1)^n(q)_n\left(1-(q^{n+2})_{\infty}\right)+2(q)_{\infty}\sum\limits_{n=0}^{\infty}\dfrac{q^n}{1+q^n}
\\
%&= 2\sum\limits_{n=0}^{\infty}(-1)^n(q)_n\left(1-(q^{n+2})_{\infty}\right)+2(q)_{\infty}\sum\limits_{n=1}^{\infty}q^n\sum\limits_{r=0}^{\infty}(-q^n)^r+(q)_{\infty}
%\\
%&= 2\sum\limits_{n=0}^{\infty}(-1)^n(q)_n\left(1-(q^{n+2})_{\infty}\right)+2(q)_{\infty}\sum\limits_{r=0}^{\infty}(-1)^r\sum\limits_{n=1}^{\infty}(q^{r+1})^n+(q)_{\infty}
%\\
&= 2\sum\limits_{n=0}^{\infty}(-1)^n(q)_n\left(1-(q^{n+2})_{\infty}\right)+2(q)_{\infty}\sum\limits_{r=0}^{\infty}\dfrac{(-1)^rq^{r+1}}{1-q^{r+1}}+(q)_{\infty}
\\
&= 2\sum\limits_{n=0}^{\infty}(-1)^n(q)_n\left(1-(q^{n+2})_{\infty}\right)+2\sum\limits_{n=0}^{\infty}(-1)^nq^{n+1}(q)_n(q^{n+2})_{\infty}+(q)_{\infty}
\\
&= 2\sum\limits_{n=0}^{\infty}(-1)^n\left((q)_n-(q)_{\infty}\right)+(q)_{\infty}
\\
&= 2(q)_{\infty}\sum\limits_{n=1}^{\infty}\dfrac{q^n}{(q)_n(1+q^n)}+(q)_{\infty}
\\
&\quad\text{(using Equation \eqn{rajat} with $(c, t) \mapsto (-1, q)$)}
\\
&= R(q), \mbox{and}
\end{align*}
\begin{align*}
B(q)&= \dfrac{\partial}{\partial z}\Big\vert_{z=1}\sum\limits_{n=0}^{\infty}\dfrac{z^nq^{\frac{n(n+1)}{2}}}{(-q)_n}
\\
&= \dfrac{\partial}{\partial z}\Big\vert_{z=1}\lim\limits_{\tau\rightarrow 0}\sum\limits_{n=0}^{\infty}\dfrac{(-q/\tau)_nz^n\tau^n}{(-q)_n}
\\
&= \dfrac{\partial}{\partial z}\Big\vert_{z=1}\lim\limits_{\tau\rightarrow 0}{}_2\phi_1\left( \begin{gathered} -q/\tau, q \\ -q \end{gathered}
;\,q, z \tau \right)
\\
&= \dfrac{\partial}{\partial z}\Big\vert_{z=1}\lim\limits_{\tau\rightarrow 0}\dfrac{(q)_{\infty}(-zq)_{\infty}}{(-q)_{\infty}(z\tau)_{\infty}}{}
_2\phi_1\left( \begin{gathered} -1, z\tau \\ -z\tau \end{gathered}
;\,q, q \right)
\\
&\quad\text{(using Equation \eqn{heine} with $(a, b, c, z) \mapsto (-q/\tau, q, -q, z\tau)$ in the first double sum)}
\\
&= \dfrac{\partial}{\partial z}\Big\vert_{z=1}\dfrac{(q)_{\infty}}{(-q)_{\infty}}\sum\limits_{n=0}^{\infty}\dfrac{(-1)_n(-zq^{n+1})_{\infty}q^n}{(q)_n}
\\
&= \dfrac{(q)_{\infty}}{(-q)_{\infty}}\sum\limits_{n=0}^{\infty}\dfrac{(-1)_n(-q^{n+1})_{\infty}q^n}{(q)_n}\sum\limits_{i=n+1}^{\infty}\dfrac{q^i}{1+q^i}
\\
&= 2(q)_{\infty}\sum\limits_{n=0}^{\infty}\dfrac{q^n}{(q)_n(1+q^n)}\sum\limits_{i=n+1}^{\infty}\dfrac{q^i}{1+q^i}.
\end{align*}
This gives us the desired result.\qed

\subsection{Proof of Theorem \thm{mexsum2}}
The proof follows from standard combinatorial arguments and using $q$-binomial theorem and Heine's ${}_2\phi_1$ transformation stated at the beginning of the section. Clearly, we have
\begin{align*}
\overline{M}_o(z,q) &= \sum\limits_{n=1}^{\infty}\dfrac{z^nq^{1+3+\cdots+(2n-3)}\prod\limits_{m=n+1}^{\infty}(1+q^{2m-1})}{\prod\limits_{\substack{m=1 \\ m \neq n}}^{\infty}(1-q^{2m-1})}
\\
&= \frac{(-q;q^2)_{\infty}}{(q;q^2)_{\infty}}\sum\limits_{n=1}^{\infty}\dfrac{z^nq^{(n-1)^2}(1-q^{2n-1})}{(-q;q^2)_n}
\\
&= \overline{P}_\text{o}(q)\left[\sum\limits_{n=1}^{\infty}\dfrac{z^nq^{(n-1)^2}}{(-q;q^2)_n}-\sum\limits_{n=1}^{\infty}\dfrac{z^nq^{n^2}}{(-q;q^2)_n}\right].
\end{align*}
Thus
\begin{align*}
\sigma \overline{M}_\text{o}(q) &= \dfrac{\partial}{\partial z}\Big\vert_{z=1} \overline{M}_o(z,q)
\\
&= \overline{P}_\text{o}(q)\left[C(q)-D(q)\right],
\end{align*}
where
\begin{align*}
C(q) &= \dfrac{\partial}{\partial z}\Big\vert_{z=1}\sum\limits_{n=1}^{\infty}\dfrac{z^nq^{(n-1)^2}}{(-q;q^2)_n}
\\
&= \dfrac{\partial}{\partial z}\Big\vert_{z=1}z\sum\limits_{n=0}^{\infty}\dfrac{z^nq^{n^2}}{(-q;q^2)_{n+1}}
\\
&= \sum\limits_{n=0}^{\infty}\dfrac{q^{n^2}}{(-q;q^2)_{n+1}}+z\dfrac{\partial}{\partial z}\Big\vert_{z=1}\sum\limits_{n=0}^{\infty}\dfrac{z^nq^{n^2}}{(-q;q^2)_{n+1}}
\\
&= 1+z\dfrac{\partial}{\partial z}\Big\vert_{z=1}\sum\limits_{n=0}^{\infty}\dfrac{z^nq^{n^2}}{(-q;q^2)_{n+1}}
\\
&\quad\text{(using Equation \eqn{gasrah} with $(z, q) \mapsto (q, q^2)$ in the first sum)}
\\
&= 1+z\dfrac{\partial}{\partial z}\Big\vert_{z=1}\sum\limits_{n=0}^{\infty}\dfrac{z^nq^{n^2}(q^2;q^2)_n}{(q^2;q^2)_n(-q;q^2)_{n+1}}
\\
&= 1+z\dfrac{\partial}{\partial z}\Big\vert_{z=1}\dfrac{(q^2;q^2)_{\infty}}{(-q;q^2)_{\infty}}\sum\limits_{n=0}^{\infty}\dfrac{z^nq^{n^2}(-q^{2n+3};q^2)_{\infty}}{(q^2;q^2)_n(q^{2n+2};q^2)_{\infty}}
\\
&= 1+z\dfrac{\partial}{\partial z}\Big\vert_{z=1}\dfrac{(q^2;q^2)_{\infty}}{(-q;q^2)_{\infty}}\sum\limits_{n=0}^{\infty}\dfrac{z^nq^{n^2}}{(q^2;q^2)_n}\sum\limits_{m=0}^{\infty}\dfrac{(-q;q^2)_mq^{(2n+2)m}}{(q^2;q^2)_m}
\\
&\quad\text{(using Equation \eqn{qbinom} with $(a, z, q) \mapsto (-q, q^{2n+2}, q^2)$)}
\\
&= 1+z\dfrac{\partial}{\partial z}\Big\vert_{z=1}\dfrac{(q^2;q^2)_{\infty}}{(-q;q^2)_{\infty}}\sum\limits_{m=0}^{\infty}\dfrac{(-q;q^2)_mq^{2m}}{(q^2;q^2)_m}\sum\limits_{n=0}^{\infty}\dfrac{z^nq^{n^2+2mn}}{(q^2;q^2)_n}
\\
&= 1+z\dfrac{\partial}{\partial z}\Big\vert_{z=1}\dfrac{(q^2;q^2)_{\infty}}{(-q;q^2)_{\infty}}\sum\limits_{m=0}^{\infty}\dfrac{(-q;q^2)_mq^{2m}(-zq^{2m+1};q^2)_{\infty}}{(q^2;q^2)_m}
\\
&\quad\text{(using \cite[Equation (2.2.6)]{An98} with $(z, q) \mapsto (zq^{2m+1}, q^2)$)}
\\
&= 1+z\dfrac{\partial}{\partial z}\Big\vert_{z=1}\dfrac{(q^2;q^2)_{\infty}(-zq;q^2)_{\infty}}{(-q;q^2)_{\infty}}\sum\limits_{m=0}^{\infty}\dfrac{(-q;q^2)_mq^{2m}}{(q^2;q^2)_m(-zq;q^2)_m}
\\
&= 1+z\dfrac{\partial}{\partial z}\Big\vert_{z=1}\dfrac{(q^2;q^2)_{\infty}(-zq;q^2)_{\infty}}{(-q;q^2)_{\infty}}{}
_2\phi_1\left( \begin{gathered} -q, 0 \\ -zq \end{gathered}
;\,q^2, q^2 \right)
\\
&= 1+z\dfrac{\partial}{\partial z}\Big\vert_{z=1}{}
_2\phi_1\left( \begin{gathered} z, q^2 \\ 0 \end{gathered}
;\,q^2, -q \right)
\\
&\quad\text{(using Equation \eqn{heine} with $(a, b, c, z, q) \mapsto (-q, 0, -zq, q^2, q^2)$)}
\\
&= 1-\sum\limits_{n=1}^{\infty}(-1)^n(q^2;q^2)_{n-1}q^n.
\end{align*}
Then,
\begin{align*}
C(-q)&=1-\sum\limits_{n=1}^{\infty}(q^2;q^2)_{n-1}q^n
\\
&= 1-q\lim_{c\rightarrow 0}{}
_2\phi_1\left( \begin{gathered} q^2,q^2 \\ c \end{gathered}
;\,q^2, q \right)
\\
&= 1-q\lim_{c\rightarrow 0}\dfrac{(c/q^2;q^2)_{\infty}(q^3;q^2)_{\infty}}{(c;q^2)_{\infty}(q;q^2)_{\infty}}{}
_2\phi_1\left( \begin{gathered} q^5/c,q^2 \\ q^3 \end{gathered}
;\,q^2, c/q^2 \right)
\\
&\quad\text{(using Equation \eqn{heine} with $(a, b, c, z, q) \mapsto (q^2, q^2, c, q, q^2)$)}
\\
&= 1-\sum\limits_{n=0}^{\infty}\dfrac{(-1)^nq^{(n+1)^2}}{(q;q^2)_{n+1}}
\\
&= 1+F(q),
\end{align*}
which gives $$C(q) = 1+F(-q),$$ and
\begin{align*}
D(q)&= \dfrac{\partial}{\partial z}\Big\vert_{z=1}\sum\limits_{n=1}^{\infty}\dfrac{z^nq^{n^2}}{(-q;q^2)_n}
\\
&= \dfrac{\partial}{\partial z}\Big\vert_{z=1}z\sum\limits_{n=0}^{\infty}\dfrac{z^nq^{(n+1)^2}}{(-q;q^2)_{n+1}}
\\
&= \sum\limits_{n=0}^{\infty}\dfrac{q^{(n+1)^2}}{(-q;q^2)_{n+1}}+z\dfrac{\partial}{\partial z}\Big\vert_{z=1}\sum\limits_{n=0}^{\infty}\dfrac{z^nq^{(n+1)^2}}{(-q;q^2)_{n+1}}
\\
&= F(-q)+z\dfrac{\partial}{\partial z}\Big\vert_{z=1}\sum\limits_{n=0}^{\infty}\dfrac{z^nq^{(n+1)^2}}{(-q;q^2)_{n+1}}
\\
&= F(-q)+z\dfrac{\partial}{\partial z}\Big\vert_{z=1}\sum\limits_{n=0}^{\infty}\dfrac{z^nq^{(n+1)^2}(q^2;q^2)_n}{(q^2;q^2)_n(-q;q^2)_{n+1}}
\\
&= F(-q)+z\dfrac{\partial}{\partial z}\Big\vert_{z=1}\dfrac{(q^2;q^2)_{\infty}}{(-q;q^2)_{\infty}}\sum\limits_{n=0}^{\infty}\dfrac{z^nq^{(n+1)^2}(-q^{2n+3};q^2)_{\infty}}{(q^2;q^2)_n(q^{2n+2};q^2)_{\infty}}
\\
&= F(-q)+z\dfrac{\partial}{\partial z}\Big\vert_{z=1}\dfrac{q(q^2;q^2)_{\infty}}{(-q;q^2)_{\infty}}\sum\limits_{n=0}^{\infty}\dfrac{z^nq^{n^2+2n}}{(q^2;q^2)_n}\sum\limits_{m=0}^{\infty}\dfrac{(-q;q^2)_mq^{(2n+2)m}}{(q^2;q^2)_m}
\\
&\quad\text{(using Equation \eqn{qbinom} with $(a, z, q) \mapsto (-q, q^{2n+2}, q^2)$)}
\\
&= F(-q)+z\dfrac{\partial}{\partial z}\Big\vert_{z=1}\dfrac{q(q^2;q^2)_{\infty}}{(-q;q^2)_{\infty}}\sum\limits_{m=0}^{\infty}\dfrac{(-q;q^2)_mq^{2m}}{(q^2;q^2)_m}\sum\limits_{n=0}^{\infty}\dfrac{z^nq^{n^2+2mn+2n}}{(q^2;q^2)_n}
\\
&= F(-q)+z\dfrac{\partial}{\partial z}\Big\vert_{z=1}\dfrac{q(q^2;q^2)_{\infty}}{(-q;q^2)_{\infty}}\sum\limits_{m=0}^{\infty}\dfrac{(-q;q^2)_mq^{2m}(-zq^{2m+3};q^2)_{\infty}}{(q^2;q^2)_m}
\\
&\quad\text{(using \cite[Equation (2.2.6)]{An98} with $(z, q) \mapsto (zq^{2m+3}, q^2)$)}
\\
&= F(-q)+z\dfrac{\partial}{\partial z}\Big\vert_{z=1}\dfrac{q(q^2;q^2)_{\infty}(-zq^3;q^2)_{\infty}}{(-q;q^2)_{\infty}}\sum\limits_{m=0}^{\infty}\dfrac{(-q;q^2)_mq^{2m}}{(q^2;q^2)_m(-zq^3;q^2)_m}
\\
&= F(-q)+z\dfrac{\partial}{\partial z}\Big\vert_{z=1}\dfrac{q(q^2;q^2)_{\infty}(-zq^3;q^2)_{\infty}}{(-q;q^2)_{\infty}}{}
_2\phi_1\left( \begin{gathered} -q,0 \\ -zq^3 \end{gathered}
;\,q^2, q^2 \right)
\\
&= F(-q)+z\dfrac{\partial}{\partial z}\Big\vert_{z=1}q{}\,
_2\phi_1\left( \begin{gathered} zq^2,q^2 \\ 0 \end{gathered}
;\,q^2, -q \right)
\\
&\quad\text{(using Equation \eqn{heine} with $(a, b, c, z, q) \mapsto (-q, 0, -zq^3, q^2, q^2)$)}
\\
&= F(-q)-q\sum\limits_{n=1}^{\infty}(-1)^n(q^2;q^2)_nq^n\sum\limits_{i=1}^{n}\dfrac{q^{2i}}{1-q^{2i}}.
\end{align*}
This gives us the desired result.\qed

\section{Concluding remarks}
We believe that there are more interesting classes of minimal excludants for overpartitions to be discovered. With suitable variations of the definitions here, it could be possible to find further relations to other mock theta functions and interesting $q$-series. During the course of our study, we have realized that not all definitions lead to interesting results. However, the results here suggest that a systematic complete study is warranted.

\section{Acknowledgments}
We express our gratitude to George Andrews for his kind interest and suggestions during the course of this project. We are also grateful to Jeremy Lovejoy for previewing a preliminary draft of this paper and his helpful comments and suggestions. We also thank the referee for several important suggestions.

\bigskip
\hrule
\bigskip

\end{document}